\def \CC {\mathbb C}
\def \QQ {\mathbb Q}
\def \RR {\mathbb R}
\def \ZZ {\mathbb Z}

\def \epsilon{\varepsilon}

\def \fine {{\hfill \qedsymbol}}

\def \si {\sigma}

\documentclass[12pt,reqno]{amsart}

\usepackage{amsmath,amssymb}

\begin{document}


\title[]{A rigidity theorem for translates of uniformly convergent Dirichlet series}

\author[]{A.PERELLI  \lowercase{and} M.Righetti }
\date{}
\maketitle

\smallskip
{\bf Abstract.} It is well known that the Riemann zeta function, as well as several other $L$-functions, is universal in the strip $1/2<\si<1$; this is certainly not true for $\sigma>1$. Answering a question of Bombieri and Ghosh, we give a simple characterization of the analytic functions approximable by translates of $L$-functions in the half-plane of absolute convergence. Actually, this is a special case of a general rigidity theorem for translates of Dirichlet series in the half-plane of uniform convergence. Our results are closely related to Bohr's equivalence theorem.



\vskip1cm
\section{Introduction}

\smallskip
In 1975, Voronin \cite{Vor/1975} discovered the following universality property of the Riemann zeta function $\zeta(s)$. Let $f(s)$ be holomorphic and non-vanishing on a closed disk $K$ inside the strip $1/2<\si<1$, and let $\epsilon>0$; then 
\[
\liminf_{T\rightarrow\infty} \frac{1}{2T}|\{\tau\in[-T,T]: \max_{s\in K}|\zeta(s+i\tau)-f(s)|<\varepsilon\}|>0.
\]
Voronin's universality theorem has been extended in several directions, in particular involving other $L$-functions in place of $\zeta(s)$, other compact sets in place of disks, and vectors of $L$-functions in place of a single $L$-function; see the survey by Matsumoto \cite{Mat/2015} and Chapter VII of Karatsuba-Voronin \cite{Ka-Vo/1992}. On the other hand, it is well known that every Dirichlet series $F(s)$ is Bohr almost periodic and bounded on any vertical strip whose closure lies inside the half-plane $\si>\si_u(F)$ of uniform convergence, hence $F(s)$ cannot be universal in the above sense for $\si>\si_u(F)$; in particular, $\zeta(s)$ is not universal for $\sigma>1$.

\medskip
In connection with their investigations on the zeros of Davenport-Heilbronn-type functions in the half-plane of absolute convergence, Bombieri-Ghosh \cite[p.230]{Bo-Gh/2011} asked for a simple characterization of the class of analytic vector functions approximable by translates of a vector of $L$-functions in the domain of absolute convergence. Here we answer this question in a rather general framework; it turns out that the answer is closely related to Bohr's theory of equivalent Dirichlet series, see Bohr \cite{Boh/1918} and Chapter 8 of Apostol \cite{Apo/1976}.

\medskip
We recall that a general Dirichlet series (D-series for short) is of the form
\begin{equation}
\label{1}
F(s) = \sum_{n=1}^\infty a(n) e^{-\lambda_n s}
\end{equation}
with coefficients $a(n)\in\CC$ and a strictly increasing sequence of real exponents $\Lambda=(\lambda_n)$ satisfying $\lambda_n\to\infty$. Clearly, the case $\lambda_n=\log n$ recovers the ordinary D-series. According to Bohr, a (possibly finite) sequence of real numbers $B=(\beta_\ell)$ is a basis of $\Lambda$ if it satisfies the following conditions: the elements of $B$ are $\QQ$-linearly independent, every $\lambda_n$ is a $\QQ$-linear combination of elements of $B$ and, viceversa, every $\beta_\ell$ is a $\QQ$-linear combination of elements of $\Lambda$. This can be expressed in matrix notation by considering $\Lambda$ and $B$ as column vectors, and writing $\Lambda = RB$ and $B=T\Lambda$ for some (infinite) Bohr matrices $R$ and $T$, whose row entries are rational and almost always 0; clearly, $R$ is uniquely determined by $\Lambda$ and $B$. Moreover, two general D-series, say $F(s)$ as in \eqref{1} and $G(s)$ with coefficients $b(n)$ and the same exponents $\Lambda$, are equivalent if there exist a basis $B$ of $\Lambda$ and a real column vector $Y=(y_\ell)$ such that
\begin{equation}
\label{2}
b(n) = a(n) e^{i(RY)_n},
\end{equation}
where $R$ is the above Bohr matrix. In the case of ordinary D-series with coefficients $a(n)$ and $b(n)$, equivalence reduces to the existence of a completely multiplicative function $\rho(n)$ such that $b(n)=a(n)\rho(n)$ for all $n\geq1$, and $|\rho(n)|=1$ whenever $a(p)\neq0$ and $p$ is a prime divisor of $n$. We refer to Chapter 8 of \cite{Apo/1976} for an introduction to Bohr's theory.

\medskip
We extend the above notion of equivalence to vectors $(F_1(s),\dots,F_N(s))$ of D-series in the following way. Let $N\geq 1$ and $F_j(s)$, $G_j(s)$, $j=1,\dots,N$, be as in \eqref{1} with coefficients $a_j(n)$ and $b_j(n)$, respectively, and the same exponents $\Lambda$. We say that $(F_1(s),\dots,F_N(s))$ and $(G_1(s),\dots,G_N(s))$ are {\it vector-equivalent} if there exist a basis $B$ of $\Lambda$ and a real vector $Y=(y_\ell)$ such that for $j=1,\dots,N$ we have
\begin{equation}
\label{3}
b_j(n) = a_j(n) e^{i(RY)_n},
\end{equation}
$R$ being as above. We stress that in \eqref{3} we require the same vector $Y$ for every $j$, hence $F_j(s)$ and $G_j(s)$ are equivalent via the same twist by $e^{i(RY)_n}$. Note that for $N=1$, vector-equivalence reduces to the standard Bohr equivalence. We also point out that we assume all the $F_j(s)$ to have the same exponents $\Lambda$ just for convenience, since otherwise we may take as $\Lambda$ the union of the exponents $\Lambda_j$ and express all the $F_j(s)$'s in terms of $\Lambda$. Moreover, as in Righetti \cite{Rig/2017}, we say that a D-series $F(s)$ as in \eqref{1}, or a sequence of exponents $\Lambda$, has an {\it integral basis} if there exists a basis $B$ of $\Lambda$ such that the associated Bohr matrix $R$ has integer entries. Such a basis $B$ is called an {\it integral basis} of $F(s)$, or of $\Lambda$. Clearly, $\Lambda=(\log n)$ has the integral basis $B=(\log p)$, so the important class of ordinary D-series falls in this case.

\medskip
Vectors of D-series with an integral basis provide a general framework where the above mentioned problem by Bombieri and Ghosh can be settled in the following sharp form. Let $N\geq 1$ and, for $j=1,\dots,N$, let $F_j(s)$ be general D-series with coefficients $a_j(n)$ and the same exponents $\Lambda$, with an integral basis and with finite $\sigma_u(F_j)$. Further, let $K_j$ be compact sets inside the half-planes $\si>\si_u(F_j)$ containing at least one accumulation point, and let $f_j(s)$ be holomorphic on $K_j$.

\medskip
{\bf Theorem 1.} {\sl Under the above assumptions, the following assertions are equivalent.

(i) For every $\epsilon>0$ there exists $\tau\in\RR$ such that
\[
\max_{j=1,\dots,N}\max_{s\in K_j} |F_j(s+i\tau)-f_j(s)|<\epsilon;
\]

(ii) $f_1(s),\dots,f_N(s)$ are general Dirichlet series with exponents $\Lambda$, and $(f_1(s),\dots,f_N(s))$ is vector-equivalent to $(F_1(s),\dots,F_N(s))$;

(iii) for every $\epsilon>0$ we have
\[
\liminf_{T\to\infty} \frac{1}{2T} |\{\tau\in[-T,T]: \max_{j=1,\dots,N}\max_{s\in K_j}|F_j(s+i\tau) - f_j(s)| <\epsilon\}| >0;
\]

(iv) $f_j(s)$ has analytic continuation to $\sigma>\sigma_u(F_j)$ and there exists a sequence $\tau_k$ such that $F_j(s+i\tau_k)$ converges uniformly to $f_j(s)$ on every closed vertical strip in $\sigma>\sigma_u(F_j)$, $j=1,\dots,N$.}

\medskip
{\bf Corollary.} {\sl Theorem $1$ holds for ordinary Dirichlet series.}

\medskip
Our result may therefore be regarded as a general {\sl rigidity} theorem for translates of D-series in the half-plane of uniform convergence, and represents the counterpart of the universality theorems for $L$-functions in the critical strip. Indeed, Theorem 1 gives a complete characterization of the analytic functions $f_j(s)$, called the {\it target functions}, approximable by such translates as in (i), and the target functions are quite special. For example, thanks to Bohr's equivalence theorem (see Theorem 8.16 of \cite{Apo/1976}) and its converse for D-series with an integral basis (see Righetti \cite{Rig/2017}), the functions $f_j(s)$ are those assuming the same set of values of the $F_j(s)$'s on any vertical strip inside the domain of absolute convergence. Moreover, if $f_j(s)$ is a target function on a compact set $K_j$ as in Theorem 1, then by (iv) it has continuation to $\sigma>\sigma_u(F_j)$ and is a target function on any compact set in such half-plane. We further note that the role of $F_j(s)$ and $f_j(s)$ in (iv), and essentially in Theorem 1, may be interchanged.

\smallskip
Note also that comparison with universality theorems  for vectors of $L$-functions is more transparent using (iii) of Theorem 1, which embodies the effect of the Kronecker-Weyl theorem. Moreover, somehow unexpectedly, contrary to the case of such universality theorems, no independence relation among the $F_j(s)$'s is required in our result. Indeed, in the special case of vectors of orthogonal $L$-functions one obtains exactly the same result as for general D-series with an integral basis. We further remark that one cannot expect Theorem 1 to hold in a larger half-plane, at least in such a general framework, since, for example, the abscissa of uniform convergence of the Dirichlet $L$-functions with primitive character equals 1, and such $L$-functions are universal in $1/2<\sigma<1$. We refer to Kaczorowski-Perelli \cite{Ka-Pe/MathZ} for a discussion of the convergence abscissae of $L$-functions.

\medskip
The interest of Bombieri and Ghosh in the above problem was related to the expectation that the real parts $\beta$ of the zeros of linear combinations of $L$-functions are dense in the interval $(1,\sigma^*)$, where $\sigma^*$ is the supremum of the $\beta$'s. However, such expectation has been shown to be incorrect by Righetti \cite{Rig/2015}, by means of counterexamples of rather general nature. The rigidity property of the translates proved in Theorem 1, and in particular the fact that the vector $Y$ in \eqref{3} is the same for all $j$'s, may possibly provide a more conceptual explanation for the existence of ``holes'' in the distribution of such real parts. However, at present we cannot make precise this assertion.

\medskip
In the next section we add some remarks on the relevance of integral bases in Theorem 1; these remarks are summarized in Theorem 2 at the end of the paper. Here we finally note that for simplicity we stated the equivalence between (i)-(iv) above under the assumption that $\Lambda$ has an integral basis, although some of the implications  hold in full generality; this will be clear from the proof.

\bigskip
\section{Proofs and remarks}

\smallskip
We need the following result about uniformly convergent D-series, which we couldn't find in the literature.

\medskip
{\bf Lemma 1.} {\sl Equivalent general Dirichlet series have the same abscissa of uniform convergence.}

\medskip
{\it Proof.} Let $F(s)$ be as in \eqref{1}; we use the following formula for $\sigma_u(F)$ due to Kuniyeda \cite{Kun/1916}. For $x\in\RR$ let
\[
T_x = \sup_{t\in\RR} \left| \sum_{[x]\leq \lambda_n<x} a(n) e^{-\lambda_n it} \right|;
\]
then
\[
\sigma_u(F) = \limsup_{x\to\infty} \frac{\log T_x}{x}.
\]
If $G(s)$ is equivalent to $F(s)$, then its coefficients $b(n)$ are given by \eqref{2}. Hence, since for fixed $x$ only finitely many $\lambda_n$'s are involved in the definition of $T_x$, we can apply Kronecker's approximation theorem to show that for every $\epsilon>0$ there exists $\tau_x\in\RR$ such that
\[
\left| \sum_{[x]\leq \lambda_n<x} a(n) e^{-\lambda_n i(t+\tau_x)} \right| - \epsilon  \leq \left| \sum_{[x]\leq \lambda_n<x} a(n)e^{i(RY)_n} e^{-\lambda_n it} \right| \leq \left| \sum_{[x]\leq \lambda_n<x} a(n) e^{-\lambda_n i(t+\tau_x)} \right| + \epsilon.
\]
See \eqref{12} and \eqref{13} at the end of the proof of Theorem 1 for details on the argument \`a la Bohr leading to the above inequalities. But
\[
\sup_{t\in\RR} \left| \sum_{[x]\leq \lambda_n<x} a(n) e^{-\lambda_n i(t+\tau_x)} \right| =  \sup_{t\in\RR} \left| \sum_{[x]\leq \lambda_n<x} a(n) e^{-\lambda_n it} \right|,
\]
and the lemma follows. \fine

\medskip
The main step in the proof of Theorem 1 is the following lemma.

\medskip
{\bf Lemma 2.} {\sl Let $F_j(s)$, $j=1,\dots,N$, be as in Theorem $1$ and let $\tau_m$ be a sequence of real numbers. Then there exists a subsequence $\tau_{m_k}$ such that, as $k\to\infty$ and for $j=1,\dots,N$, $F_j(s+i\tau_{m_k})$ converges uniformly on any closed vertical strip inside $\si>\si_u(F_j)$ to a general Dirichlet series $G_j(s)$ with exponents $\Lambda$, and $(G_1(s),\dots,G_N(s))$ is vector-equivalent to $(F_1(s), \dots,F_N(s))$.}

\medskip
{\it Proof.} Let $B=(\beta_\ell)$ be an integral basis of the exponents $\Lambda$ of the $F_j(s)$, and let
\[
\theta_{m,\ell}=\left\{-\frac{\tau_m\beta_\ell}{2\pi}\right\},\qquad m,\ell=1,2,\ldots,
\]
where $\{x\}$ denotes the fractional part of $x$. Since $0\leq \theta_{m,\ell}<1$, by Helly's selection principle, see Lemma 1 of Section 8.12 of \cite{Apo/1976}, there exist a subsequence $m_k$ and a sequence of real numbers $\theta_\ell$ such that 
\begin{equation}
\label{4}
\lim_{k\rightarrow\infty} \theta_{m_k,\ell}=\theta_\ell
\end{equation}
for every $\ell\geq 1$. Next we define $Y=(2\pi \theta_\ell)$ and, for $j=1,\dots,N$,
\begin{equation}
\label{5}
G_j(s)=\sum_{n=1}^\infty a_j(n) e^{i(RY)_n} e^{-\lambda_ns},
\end{equation}
where $R=(r_{n,\ell})$ is the Bohr matrix such that $\Lambda=RB$. Clearly, $(G_1(s),\dots,G_N(s))$ is vector-equivalent to $(F_1(s), \dots,F_N(s))$ by definition, and now we show that every $F_j(s+i\tau_{m_k})$ converges to $G_j(s)$ uniformly over any closed vertical strip inside $\sigma>\sigma_u(F_j)$.

\smallskip
We first note that since $B$ is an integral basis of $\Lambda$ we have
\[
e^{-i\lambda_n\tau_{m_k}} = e^{2\pi i\sum_\ell r_{n,\ell}(-\frac{\tau_{m_k}\beta_\ell}{2\pi})} = e^{2\pi i (\sum_\ell r_{n,\ell}\theta_{m_k,\ell})},
\]
hence
\begin{equation}
\label{6}
e^{-i\lambda_n\tau_{m_k}} - e^{i(RY)_n} = e^{2\pi i( \sum_\ell r_{n,\ell} \theta_\ell)} \left(e^{2\pi i \sum_\ell r_{n,\ell} (\theta_{m_k,\ell}-  \theta_\ell)} - 1\right).
\end{equation}
Moreover, recalling that the row entries of $R$ are almost always 0, for every $n\geq1$ there exists $c_n\geq1$ such that
\begin{equation}
\label{7}
\left|\sum_\ell r_{n,\ell} (\theta_{m_k,\ell} - \theta_\ell) \right|  \leq  c_n \max_{\ell \, \text{with} \, r_{n,\ell}\neq 0} |\theta_{m_k,\ell}-\theta_\ell|.
\end{equation}
Let now $W_j$ be a closed vertical strip inside $\si>\si_u(F_j)$, and let $\epsilon>0$ be sufficiently small. By the uniform convergence and thanks to Lemma 1, there exists $M=M_j(\epsilon)$ such that 
\begin{equation}
\label{8}
\sup_{s\in W_j}\left(\left|\sum_{n>M} a_j(n) e^{-\lambda_n(s+i\tau_{m_k})}\right|, \left|\sum_{n>M} a_j(n) e^{i(RY)_n}  e^{-\lambda_ns}\right|\right)< \epsilon.
\end{equation}
Next, writing
\[
C=C_j(\epsilon) = \max_{n\leq M} c_n \quad \text{and} \quad H=H_j(\epsilon) = \max_{s\in W_j} \sum_{n\leq M}|a_j(n)| e^{-\lambda_n\sigma},
\]
in view of \eqref{4} there exists $\overline{k}=\overline{k}_j(\epsilon)$ such that for $k\geq \overline{k}$
\begin{equation}
\label{9}
\max_{\ell \, \text{with} \, r_{n,\ell}\neq 0} |\theta_{m_k,\ell}-\theta_\ell|<\epsilon/CH
\end{equation}
for every $n\leq M$. Hence, from \eqref{6}-\eqref{9}, for $k\geq \overline{k}$ we have that
\begin{equation}
\label{10}
\begin{split}
\sup_{s\in W_j} |&F_j(s+i\tau_{m_k}) - G_j(s)| \\
& < 2\epsilon + \max_{s\in W_j} \sum_{n\leq M} |a_j(n)| \left| e^{-i\lambda_n\tau_{m_k}} - e^{i(RY)_n} \right| e^{-\lambda_n\sigma} < 10\epsilon,
\end{split}
\end{equation}
and the lemma follows. \fine

\medskip
{\bf Proof of Theorem 1.} From (i) applied with $\epsilon=1/m$, $m=1,2,\dots$, we obtain a sequence $\tau_m$ such that $F_j(s+i\tau_m)$ converges uniformly to $f_j(s)$ over $K_j$, for $j=1,\dots,N$. Thanks to Lemma 2 there exists a subsequence $\tau_{m_k}$ such that $F_j(s+i\tau_{m_k})$ converges uniformly over $K_j$ to $G_j(s)$. Hence $f_j(s)=G_j(s)$ by the uniqueness of the limit and of the analytic continuation, and (ii) follows from the properties of the $G_j(s)$'s in Lemma 2.

\smallskip
Suppose now that the $f_j(s)$'s are as in (ii), hence their coefficients $b_j(n)$ are as in \eqref{3} with the same $Y=(y_\ell)$, and let $R=(r_{n,\ell})$ be the Bohr matrix of a basis $B=(\beta_\ell)$ of $\Lambda$. Note that {\it here we do not assume that $\Lambda$ has an integral basis and that the $K_j$'s have an accumulation point.} Given $\epsilon>0$ and $\tau\in\RR$, thanks to Lemma 1 let, as in the proof of Lemma 2, $M=M(\epsilon)>0$ be such that
\begin{equation}
\label{11}
\begin{split}
\max_{j=1,\dots,N} \max_{s\in K_j} &|F_j(s+i\tau)-f_j(s)| \\
&<2\epsilon +\max_{j=1,\dots,N} \max_{s\in K_j} \sum_{n\leq M}|a_j(n)| \left|e^{-i\lambda_n\tau} - e^{i(RY)_n} \right| e^{-\lambda_n\sigma}.
\end{split}
\end{equation}
Recalling the properties of the Bohr matrices, we express the exponents $\lambda_n$ by means of the basis $B$, write $r_{n,\ell}=a_{n,\ell}/q_{n,\ell}$ and finally denote by $Q=Q(\epsilon)$ the least common multiple of all the $q_{n,\ell}$'s, with $n\leq M$ and $\ell \geq 1$, such that $r_{n,\ell}\neq0$. We thus obtain, for $n\leq M$, that
\begin{equation}
\label{12}
e^{-i\lambda_n\tau} - e^{i(RY)_n} = e^{2\pi i \sum_\ell m_{n,\ell}(\frac{y_\ell}{2\pi Q})} \big(e^{2\pi i \sum_\ell m_{n,\ell} (-\frac{\beta_\ell\tau}{2\pi Q} -  \frac{y_\ell}{2\pi Q})} - 1\big)
\end{equation}
with certain $m_{n,\ell}\in\ZZ$. Since the $\beta_\ell$ are $\QQ$-linearly independent, by Kronecker's approximation theorem (see e.g. Chapter 8 of Chandrasekharan \cite{Cha/1968}) for every $\delta>0$ there exists $\tau\in\RR$ such that
\begin{equation}
\label{13}
\left\|-\frac{\beta_\ell\tau}{2\pi Q} - \frac{y_\ell}{2\pi Q} \right\| < \delta
\end{equation}
for all $\ell$ involved in \eqref{12} with $n\leq M$, where $\|x\|$ denotes the distance of $x$ from the nearest integer. As in Lemma 2, by an obvious choice of $\delta$ in terms of $\epsilon$, of $F_j(s)$ and $K_j$ for $j=1,\dots,N$ and of $\max_{n\leq M} \sum_\ell |m_{n,\ell}|$, from \eqref{11}-\eqref{13} we obtain that there exists $\tau\in\RR$ such that
\[
\max_{j=1,\dots,N}\max_{s\in K_j}|F_j(s+i\tau)-f_j(s)| \ll \epsilon,
\]
and (i) follows.

\smallskip
Finally, clearly (iii) implies (i), and replacing Kronecker's approximation theorem by the Kronecker-Weyl theorem (see Appendix 8 of \cite{Ka-Vo/1992} or Remark 1.1 on p.96-97 in \cite{Mat/2015}) in the above proof that (ii) implies (i), we can show that (ii) implies (iii) as well. Moreover, clearly (iv) implies (i), while (i) implies (iv) thanks to Lemma 2 exactly as in the above proof that (i) implies (ii), choosing $\tau_k=\tau_{m_k}$. The proof of Theorem 1 is now complete. \fine

\bigskip
We conclude with some remarks about the relevance of integral bases in Theorem 1. We already remarked that the D-series with an integral basis contain the ordinary D-series. A simple but interesting example of non-ordinary D-series with an integral basis is the Hurwitz zeta function
\[
\sum_{n=0}^\infty \frac{1}{(n+\alpha)^s}
\]
with a transcendental $0<\alpha<1$. Indeed, in this case the exponents $\lambda_n=\log(n+\alpha)$ are all $\QQ$-linearly independent, see Davenport-Heilbronn \cite{Da-Ha/1936}, therefore $\Lambda$ is already a basis and hence $R$ is the identity matrix.

\medskip
Even if $\Lambda$ does not have an integral basis, it is still possible to say something on the target functions $f_j(s)$ by a variant of the above arguments, although such a set may be larger in this case since we have seen that (ii) implies (i) in full generality. From now on we assume (i) as in Theorem 1, but not anymore that $\Lambda$ has an integral basis. We first note that by a variant of the first steps of Lemma 2, namely considering the double sequence
\[
\theta_{m,n}=\left\{-\frac{\tau_m\lambda_n}{2\pi}\right\},\qquad m,n=1,2,\ldots
\]
and the sequence $\theta_n$ obtained as in \eqref{4}, we are led to the D-series
\begin{equation}
\label{14}
G_j(s) = \sum_{n=1}^\infty a_j(n)e^{2\pi i\theta_n} e^{-\lambda_n s}, \qquad j=1,\dots,N,
\end{equation}
instead of those in \eqref{5}. Next, we observe that a (simpler) variant of Lemma 1 shows that $\sigma_u(G_j)=\sigma_u(F_j)$, for $j=1,\dots,N$. Indeed, for every $\epsilon>0$ there exists $k=k(x)$ such that
\[
\left| \sum_{[x]\leq \lambda_n<x} a(n) e^{-\lambda_n i(t+\tau_{m_k})} \right| - \epsilon  \leq \left| \sum_{[x]\leq \lambda_n<x} a(n)e^{2\pi i\theta_n} e^{-\lambda_n it} \right| \leq \left| \sum_{[x]\leq \lambda_n<x} a(n) e^{-\lambda_n i(t+\tau_{m_k})} \right| + \epsilon,
\]
and the assertion follows as before. Hence, by a (simpler) variant of the arguments in the second part of the proof of Lemma 2, see \eqref{8}-\eqref{10}, we obtain that $F_j(s+i\tau_{m_k})$ converges uniformly to $G_j(s)$ on any closed vertical strip inside $\si>\si_u(F_j)$, $j=1,\dots,N$. Now, having (i), it is not difficult to conclude as before that the $f_j(s)$'s coincide with the $G_j(s)$'s in \eqref{14}. In particular, $f_j(s)$ and $F_j(s)$ have the same abscissae of absolute and uniform convergence.

\smallskip
One can show that the $f_j(s)$'s have further properties; for example, denoting by $S_f(V)$ the set of values taken by $f(s)$ on $V$, we have that $S_{f_j}(V_j) \subseteq S_{F_j}(V_j)$ for any open vertical strip $V_j$ in $\si>\si_u(F_j)$, $j=1,\ldots,N$. Indeed, suppose that $v_j\in S_{f_j}(V_j)$, and that $f_j(s_j)=v_j$ for some $s_j\in V_j$; moreover, let $r_j>0$ be such that the disk $K_j=\{|s-s_j|\leq r_j\}$ is contained in $V_j$. By the above argument we know that $F_j(s+i\tau_{m_k})$ converges uniformly to $f_j(s)$ over $K_j$. If $f_j(s)$ is constant then, by \eqref{14}, $F_j(s)$ is also constant and the assertion follows trivially. Otherwise, taking $r_j$ sufficiently small we have
\[
\min_{|s-s_j|=r_j} |f_j(s)-v_j|=\eta_j>0,
\]
and certainly there exists $k$ such that
\[
\max_{|s-s_j|\leq r_j} |F_j(s+i\tau_{m_k})- f_j(s)| <\eta_j.
\]
Therefore, by an application of Rouch\'e's theorem we deduce that $F_j(s)=v_j$ has solutions for $s\in K_j$, and our assertion follows.

\smallskip
Actually, the opposite inclusion holds as well, namely $S_{F_j}(V_j) \subseteq S_{f_j}(V_j)$ for every such $V_j$. Indeed, still thanks to the above argument ensuring the uniform convergence of $F_j(s+i\tau_{m_k})$ to $f_j(s)$ over any closed vertical strip in $\sigma>\sigma_u(F_j)$, we may invert the role of $F_j(s)$ and $f_j(s)$. Therefore, for $j=1,\dots,N$, $f_j(s-i\tau_{m_k})$ converges uniformly to $F_j(s)$ on a suitable disk $K_j$ around a point $s_j$ such that $F_j(s_j)=v_j\in S_{F_j}(V_j)$, and we may conclude as before that $S_{F_j}(V_j) \subseteq S_{f_j}(V_j)$.

\medskip
Summarizing, with the above notation we have the following result.

\medskip
{\bf Theorem 2.} {\sl Under the assumptions of Theorem $1$, with $\Lambda$ not necessarily having an integral basis, suppose that (i) holds. Then the $f_j(s)$'s are general Dirichlet series with coefficients $b_j(n)$ and the same exponents $\Lambda$, and satisfy the following properties. For $j=1,\dots,N$
\[
|b_j(n)| = |a_j(n)|, \quad  \sigma_u(f_j)=\sigma_u(F_j) \quad \text{and} \quad S_{f_j}(V_j) = S_{F_j}(V_j),
\]
where $V_j$ is any open vertical strip inside $\sigma>\sigma_u(F_j)$. Moreover, (i) holds for the $f_j(s)$'s described in (ii) of Theorem $1$.}

\medskip
Similar remarks and variants, namely without assuming the existence of an integral basis, apply also to the equivalence of (i) with (iii) and (iv) in Theorem 1. However, $f_j(s)$ may not be equivalent to $F_j(s)$, as shown by the following example by Bohr \cite[pp.151--153]{Boh/1918}. Let 
\[
\lambda_n=2n-1+\frac{1}{2(2n-1)}, \qquad F(s) = \sum_{n= 1}^\infty e^{-\lambda_n s}, \qquad f(s)=-F(s).
\] 
In this case all bases $B$ of $\Lambda$ consist of a single rational number, and since the least common multiple of the denominators of the $\lambda_n$ is $\infty$, no one is an integral basis. Moreover, the Bohr matrix $R$ such that $\Lambda=RB$ reduces to an infinite column vector, hence the vectors $Y$ in \eqref{2} reduce to a single real number; thus the set of D-series equivalent to $F(s)$ consists of its vertical shifts. Further, as shown by Bohr, $f(s)$ is not equivalent to $F(s)$.  On the other hand, $f(s)$ satisfies (i) in Theorem 2 with $\tau = 2\pi \prod_{n\leq m} (2n-1)$, for any sufficiently large $m=m(\epsilon)$.

\bigskip
{\bf Acknowledgements.}  This research was partially supported by PRIN2015 {\sl Number Theory and Arithmetic Geometry}. A.P. is member of the GNAMPA group of INdAM, and M.R. was partially supported by a research scholarship of the Department of Mathematics, University of Genova.

\newpage

\ifx\undefined\bysame{poly}.
\newcommand{\bysame}{\leavevmode\hbox to3em{\hrulefill}\ ,}
\fi

\vskip1cm
Alberto Perelli, Dipartimento di Matematica, Via Dodecaneso 35, 16146 Genova, Italy. 

perelli@dima.unige.it

\medskip
Mattia Righetti, Centre de Recherches Math\'ematiques, Universit\'e de Montr\'eal, P.O. box 

6128, Centre-Ville Station, Montr\'eal, Qu\'ebec H3C 3J7, Canada.

righetti@crm.umontreal.ca

\end{document}